\RequirePackage{fix-cm}
\documentclass[smallextended]{svjour3}       
\smartqed  
\usepackage[square,numbers,sort&compress]{natbib}
\usepackage[]{geometry}
\usepackage{graphicx}
\usepackage{amsfonts}
\spnewtheorem{fact}{Fact}{\bf}{\it}
\begin{document}

\title{{Connection between
the clique number and the Lagrangian of $3$-uniform hypergraphs}
}


\author{Qingsong Tang         \and   Yuejian  Peng  \and   Xiangde Zhang \and Cheng Zhao
}


\institute{Qingsong Tang \at
              College of Sciences, Northeastern University, Shenyang, 110819, P.R.China. \at
              School of Mathematics, Jilin University, Changchun, 130012, P.R.China.\\
              \email{t\_qsong@sina.com}
              \and
              Yuejian  Peng \at
              College of Mathematics, Hunan University, Changsha 410082, P.R. China. \at This research is supported by National Natural Science Foundation of China (No. 11271116).\\
              \email{ypeng1@163.com}
               \and
           Xiangde Zhang \at
             College of Sciences, Northeastern University, Shenyang, 110819, P.R.China\\
             \email{zhangxdneu@163.com}
              \and
             Cheng Zhao(Corresponding author) \at
             Department of Mathematics and Computer Science, Indiana State University, Terre Haute, IN, 47809 USA. \at School of Mathematics, Jilin University, Changchun 130012, P.R. China.\\
             \email{cheng.zhao@indstate.edu}
}

\date{Received: date / Accepted: date}

\maketitle

\begin{abstract}
There is a remarkable connection between
the clique number and the Lagrangian of a 2-graph proved by Motzkin and Straus in 1965.
It is useful in practice if similar results hold for hypergraphs. However the obvious generalization of
Motzkin and Straus' result to hypergraphs is false. Frankl and F\"{u}redi conjectured that the $r$-uniform hypergraph with $m$ edges formed by taking the first $m$ sets in the colex ordering of ${\mathbb N}^{(r)}$ has the largest Lagrangian of all $r$-uniform hypergraphs with $m$ edges. For $r=2$, Motzkin and Straus' theorem confirms this conjecture. For $r=3$, it is shown by Talbot that this conjecture is true when $m$ is in certain ranges. In this paper, we explore  the connection between the clique number and Lagrangians for $3$-uniform hypergraphs.
As an application of this connection, we confirm that Frankl and F\"{u}redi's conjecture holds for bigger ranges of $m$ when $r$=3. We also  obtain two weaker versions of Tur\'{a}n type theorem for left-compressed $3$-uniform hypergraphs.
\keywords{Cliques of hypergraphs \and Colex ordering \and Lagrangians of hypergraphs \and  Polynomial optimization.}
\subclass{05C35 \and 05C65 \and 05D99 \and 90C27}
\end{abstract}

\section{Introduction}
\label{intro}
For a set $V$ and a positive integer $r$, let  $V^{(r)}$ be  the family of all $r$-subsets of $V$. An $r$-uniform graph or $r$-graph $G$ consists of a set $V(G)$ of vertices and a set $E(G) \subseteq V(G) ^{(r)}$ of edges. An edge $e=\{a_1, a_2, \ldots, a_r\}$ will be simply denoted by $a_1a_2 \ldots a_r$. An $r$-graph $H$ is  a {\it subgraph} of an $r$-graph $G$, denoted by $H\subseteq G$ if $V(H)\subseteq V(G)$ and $E(H)\subseteq E(G)$. Let $K^{(r)}_t$ denote the complete $r$-graph on $t$ vertices, that is the $r$-graph on $t$ vertices containing all possible edges. A complete $r$-graph on $t$ vertices is also called a clique with order $t$.  A clique is said to be maximal if there is no other clique containing it, while it is called maximum if it has maximum cardinality. The clique number of a $r$-graph $G$, denoted as $\omega(G)$, is defined as the cadinality of the  maximum clique. Let ${\mathbb N}$ be the set of all positive integers.  For an integer $n \in {\mathbb N}$, let $[n]$ denote the set $\{1, 2, 3, \ldots, n\}$. Let $[n]^{(r)}$  represent the  complete $r$-graph on the vertex set $[n]$. When $r=2$, an $r$-graph is a simple graph.  When $r\ge 3$,  an $r$-graph is often called a hypergraph.

For an $r$-graph $G:=(V,E)$, denote the $(r-1)$-neighborhood of a vertex $i \in V$ by $E_i:=\{A \in V^{(r-1)}: A \cup \{i\} \in E\}$. Similarly,  denote the $(r-2)$-neighborhood of a pair of vertices $i,j \in V$ by $E_{ij}:=\{B \in V^{(r-2)}: B \cup \{i,j\} \in E\}$.  Denote the complement of $E_i$ by $E^c_i:=\{A \in V^{(r-1)}: A \cup \{i\} \in V^{(r)} \backslash E\}$. Also, denote the complement of $E_{ij}$ by
$E^c_{ij}:=\{B \in V^{(r-2)}: B \cup \{i,j\} \in V^{(r)} \backslash E\}$.
Denote $E_{i\setminus j}:=E_i\cap E^c_j.$
An $r$-graph $G=([n],E)$ is {\it left-compressed} if $j_1j_2 \cdots j_r \in E$ implies $i_1i_2 \cdots i_r \in E$ provided $i_p \le j_p$ for every $p, 1\le p\le r$. Equivalently, an $r$-graph $G=([n],E)$ is {\it left-compressed} if $E_{j\setminus i}=\emptyset$ for any $1\le i<j\le n$.

\begin{definition} \label{definitionlagrangian}
For  an $r$-uniform graph $G$ with the vertex set $[n]$, edge set $E(G)$, and a vector \\$\vec{x}=(x_1,\ldots,x_n) \in {\mathbb R}^n$,
we associate a homogeneous polynomial in $n$ variables, denoted by  $\lambda (G,\vec{x})$  as follows:
$$\lambda (G,\vec{x})=\sum_{i_1i_2 \cdots i_r \in E(G)}x_{i_1}x_{i_2}\ldots x_{i_r}.$$
Let $S=\{\vec{x}=(x_1,x_2,\ldots ,x_n): \sum_{i=1}^{n} x_i =1, x_i
\ge 0 {\rm \ for \ } i=1,2,\ldots , n \}$.
Let $\lambda (G)$ represent the maximum
 of the above homogeneous  multilinear polynomial of degree $r$ over the standard simplex $S$. Precisely
 $$\lambda (G) = \max \{\lambda (G, \vec{x}): \vec{x} \in S \}.$$
\end{definition}
The value $x_i$ is called the {\em weight} of the vertex $i$.
A vector $\vec{x}:=(x_1, x_2, \ldots, x_n) \in {\mathbb R}^n$ is called a feasible weighting for $G$ iff
$\vec{x}\in S$. A vector $\vec{y}\in S$ is called an {\em optimal weighting} for $G$
iff $\lambda (G, \vec{y})=\lambda(G)$. We call $\lambda(G)$ the Graph-Lagrangian of hypergraph $G$, for abbreviation, the Lagrangian of $G$.

The following fact is easily implied by Definition \ref{definitionlagrangian}.
\begin{fact}\label{mono}
Let $G_1$, $G_2$ be $r$-uniform graphs and $G_1\subseteq G_2$. Then $\lambda (G_1) \le \lambda (G_2).$
\end{fact}
The maximum clique problem is a classical problem in combinatorial optimization which has
important applications in various domains. In \cite{MS}, Motzkin and Straus established a remarkable connection between the clique number and the Lagrangian of a graph.

\begin{theorem} (\cite{MS}) \label{MStheo}
If $G$ is a 2-graph with clique number $t$ then
$\lambda(G)=\lambda(K^{(2)}_t)={1 \over 2}(1 - {1 \over t})$.
\end{theorem}

The obvious generalization of Motzkin and Straus' result to hypergraphs is false because there are many examples of hypergraphs that do not achieve their Lagrangian on any proper subhypergraph. Lagrangians of hypergraphs has been proved to be a useful tool,  for example, it is useful to hypergraph  extremal problems. Applications of Lagrangian method can be found in \cite{FF, FR84, keevash, mubayi06, sidorenko89}. In most applications, an upper bound is needed. Frankl and F\"uredi \cite{FF} asked the following question. Given $r \ge 3$ and $m \in {\mathbb N}$ how large can the Lagrangian of an $r$-graph with $m$ edges be?
For distinct $A, B \in {\mathbb N}^{(r)}$ we say that $A$ is less than $B$ in the {\em colex ordering} if $max(A \triangle B) \in B$, where $A \triangle B=(A \setminus B)\cup (B \setminus A)$. For example we have $246 < 156$ in ${\mathbb N}^{(3)}$ since $max(\{2,4,6\} \triangle \{1,5,6\}) \in \{1,5,6\}$. In colex ordering, $123<124<134<234<125<135<235<145<245<345<126<136<236<146<246<346<156<256<356<456<127<\cdots .$ Note that the first $t \choose r$ $r$-tuples in the colex ordering of ${\mathbb N}^{(r)}$ are the edges of $[t]^{(r)}$. The following conjecture of Frankl and F\"uredi (if it is true) proposes a  solution to the question mentioned above.

\begin{conjecture} (\cite{FF})\label{conjecture} The $r$-graph with $m$ edges formed by taking the first $m$ sets in the colex ordering of ${\mathbb N}^{(r)}$ has the largest Lagrangian of all $r$-graphs with  $m$ edges. In particular, the $r$-graph with $t \choose r$ edges and the largest Lagrangian is $[t]^{(r)}$.
\end{conjecture}

This conjecture is true when $r=2$ by Theorem \ref{MStheo}. For the case $r=3$, Talbot in \cite{T} proved the following.

\begin{theorem} (\cite{T}) \label{Tal} Let $m$ and $t$ be integers satisfying
${t-1 \choose 3} \le m \le {t-1 \choose 3} + {t-2 \choose 2} - (t-1).$
Then Conjecture \ref{conjecture} is true for $r=3$ and this value of $m$.
Conjecture \ref{conjecture} is also true for $r=3$ and $m= {t \choose 3}-1$ or $m={t \choose 3} -2$.
\end{theorem}
Further evidence that supports Conjecture \ref{conjecture} can be found in \cite{TPWP, TPZZ1}. In particular, Conjecture \ref{conjecture} is  true for $r=3$ and ${t \choose 3}-6 \le m \le {t \choose 3}$  (see \cite{TPWP, TPZZ1}).


 Although the obvious generalization of Motzkin and Straus' result to hypergraphs is false,  we attempt to explore the relationship between the Lagrangian of a hypergraph and its cliques number for hypergraphs when the number of edges is in certain ranges. In \cite{PZ}, it is conjectured that the following Motzkin and Straus type results are true for hypergraphs.

\begin{conjecture} \label{conjecture1}
Let $t$, $m$, and $r\ge 3$ be positive integers satisfying ${t-1 \choose r} \le m \le {t-1 \choose r} + {t-2 \choose r-1}$.
Let $G$ be an $r$-graph with $m$ edges and $G$ contain a clique of order  $t-1$. Then $\lambda(G)=\lambda([t-1]^{(r)})$.
\end{conjecture}

The upper bound ${t-1 \choose r} + {t-2 \choose r-1}$ in this conjecture is the best possible. When $m ={t-1 \choose r}+{t-2 \choose r-1}+1$, let  $C_{r,m}$ be the $r$-graph with the  vertex set $[t]$ and the edge set $[t-1]^{(r)}\cup\{i_1\cdots i_{r-1}t: i_1\cdots i_{r-1}\in [t-2]^{(r-1)}\}\cup\{1\cdots (r-2)(t-1)t\}$. Take a legal weighting $\vec{x}:=(x_1,\ldots,x_t)$, where $x_1=x_2=\cdots=x_{t-2}={1 \over t-1}$ and $x_{t-1}=x_{t}={1 \over 2(t-1)}$.
Then $\lambda(C_{r,m})\ge \lambda(C_{r,m},\vec{x} )> \lambda([t-1]^{(r)})$.

\begin{conjecture} \label{conjecture2} Let $G$ be an $r$-graph with $m$ edges without containing a clique of size  $t-1$, where ${t-1 \choose r} \le m \le {t-1 \choose r} + {t-2 \choose r-1}$. Then $\lambda(G) < \lambda([t-1]^{(r)})$.
\end{conjecture}

Let $C_{r,m}$ denote the $r$-graph with $m$ edges formed by taking the first $m$ sets in the colex ordering of ${\mathbb N}^{(r)}$.
The following result was given in  \cite{T}.

\begin{lemma} \cite{T} \label{LemmaTal7}
For any integers $m, t,$ and $r$ satisfying ${t-1 \choose r} \le m \le {t-1 \choose r} + {t-2 \choose r-1}$,
we have $\lambda(C_{r,m}) = \lambda([t-1]^{(r)})$.
\end{lemma}

\begin{remark}  

 Conjectures \ref{conjecture1} and \ref{conjecture2} refine   Conjecture \ref{conjecture} when ${t-1 \choose r} \le m \le {t-1 \choose r} + {t-2 \choose r-1}$. If Conjectures \ref{conjecture1} and \ref{conjecture2} are true, then Conjecture \ref{conjecture} is true for this range of $m$.
\end{remark}

In \cite{PZ}, we showed that Conjecture \ref{conjecture1} holds when $r=3$ as in the following Theorem.

\begin{theorem} (\cite{PZ})\label{TheoremPZ18}  Let $m$ and $t$ be positive integers satisfying ${t-1 \choose 3} \le m \le {t-1 \choose 3} + {t-2 \choose 2}$. Let $G$ be a $3$-graph with $m$ edges and  contain a clique of order  $t-1$. Then $\lambda(G) = \lambda([t-1]^{(3)})$.
\end{theorem}

In this paper, we will show the following.

\begin{theorem} \label{theorem 1} Let $m$ and $t$ be integers satisfying
${t-1 \choose 3} \le m \le {t-1 \choose 3} + {t-2 \choose 2} - \frac{1}{2}(t-1).$
Let $G$ be a $3$-graph with $m$ edges without containing a clique order of $t-1$, then $\lambda(G)<\lambda([t-1]^{(3)}).$
\end{theorem}

Combing Theorems \ref{TheoremPZ18} and \ref{theorem 1}, we have the follow result on Conjecture \ref{conjecture}.
\begin{corollary} \label{improvement} Let $m$ and $t$ be integers satisfying
${t-1 \choose 3} \le m \le {t-1 \choose 3} + {t-2 \choose 2} - \frac{1}{2}(t-1).$
Then Conjecture \ref{conjecture} is true for $r=3$ and this value of $m$.
\end{corollary}

Note that Theorem \ref{theorem 1} supports Conjecture \ref{conjecture2} and Corollary \ref{improvement} improves Thoerem \ref{Tal}.

The rest of the paper is organized as follows. In section \ref{sec:2}, we prove Theorem \ref{theorem 1}.  In section \ref{sec:3}, we explore the connection between the clique number and the Lagrangians of some left-compressed 3-graphs. As an application, we obtain two weaker versions of Tu\'{a}n type result.  
First we give some useful results.
\section{Useful Results}
\label{sec:1}
We will impose one additional condition on any optimal weighting ${\vec x}=(x_1, x_2, \ldots, x_n)$ for an $r$-graph $G$:
\begin{eqnarray}
 &&|\{i : x_i > 0 \}|{\rm \ is \ minimal, i.e. \ if}  \ \vec y {\rm \ is \ a \ legal \ weighting \ for \ } G  {\rm \ satisfying }\nonumber \\
 &&|\{i : y_i > 0 \}| < |\{i : x_i > 0 \}|,  {\rm \  then \ } \lambda (G, {\vec y}) < \lambda(G) \label{conditionb}.
\end{eqnarray}

When the theory of Lagrange multipliers is applied to find the optimum of $\lambda(G, {\vec x})$, subject to $\sum_{i=1}^n x_i =1$, notice that $\lambda (E_i, {\vec x})$ corresponds to the partial derivative of  $\lambda(G, \vec x)$ with respect to $x_i$. The following lemma gives some necessary conditions of an optimal weighting for $G$.

\begin{lemma} (\cite{FR84}) \label{LemmaTal5} Let $G:=(V,E)$ be an $r$-graph on the vertex set $[n]$ and ${\vec x}=(x_1, x_2, \ldots, x_n)$ be an optimal  weighting for $G$ with $k$  ($\le n$) non-zero weights $x_1$, $x_2$, $\cdots$, $x_k$  satisfying condition (\ref{conditionb}). Then for every $\{i, j\} \in [k]^{(2)}$, (a) $\lambda (E_i, {\vec x})=\lambda (E_j, \vec{x})=r\lambda(G)$, (b) there is an edge in $E$ containing both $i$ and $j$.
\end{lemma}

\begin{remark}\label{r1} (a) In Lemma \ref{LemmaTal5}, part(a) implies that
$$x_j\lambda(E_{ij}, {\vec x})+\lambda (E_{i\setminus j}, {\vec x})=x_i\lambda(E_{ij}, {\vec x})+\lambda (E_{j\setminus i}, {\vec x}).$$
In particular, if $G$ is left-compressed, then
$$(x_i-x_j)\lambda(E_{ij}, {\vec x})=\lambda (E_{i\setminus j}, {\vec x})$$
for any $i, j$ satisfying $1\le i<j\le k$ since $E_{j\setminus i}=\emptyset$.

(b) If  $G$ is left-compressed, then for any $i, j$ satisfying $1\le i<j\le k$,
\begin{equation}\label{enbhd}
x_i-x_j={\lambda (E_{i\setminus j}, {\vec x}) \over \lambda(E_{ij}, {\vec x})}
\end{equation}
holds.  If  $G$ is left-compressed and  $E_{i\setminus j}=\emptyset$ for $i, j$ satisfying $1\le i<j\le k$, then $x_i=x_j$.

(c) By (\ref{enbhd}), if  $G$ is left-compressed, then an optimal legal weighting  ${\vec x}=(x_1, x_2, \ldots, x_n)$ for $G$  must satisfy
\begin{equation}\label{conditiona}
x_1 \ge x_2 \ge \ldots \ge x_n \ge 0.
\end{equation}
\end{remark}
The following lemma implies that we only need to consider left-compressed $r$-graphs when Conjecture \ref{conjecture} is explored.
\begin{lemma} (\cite{T}) \label{LemmaTal8}
Let $m,t$ be positive integers satisfying  $m\leq{t \choose r}-1$, then there exists a left-compressed $r$-graph $G$ with $m$ edges such that $\lambda(G)=\lambda_{m}^{r}.$
\end{lemma}
\section{Proof of Theorem \ref{theorem 1}}
\label{sec:2}
The following lemma showed in \cite{PZZ} implies that we only need to consider left-compressed  $3$-graphs $G$ on $t$ vertices  to verify Conjecture \ref{conjecture2} for $r=3$.
Denote

$\lambda_{(m,t)}^{r-}:=\max \{ \lambda(G):$
$ G$ is an  $r$ - graph  with $m$   edges  and does not  contain a clique of size    $t  \}$ .
\begin{lemma} \label{lemma1}\cite{PZZ} Let $m$ and $t$  be  positive integers satisfying ${t-1 \choose 3} \le m \le {t-1 \choose 3} + {t-2 \choose 2}$.
Then there exists a left-compressed $3$-graph $G$ on the vertex set $[t]$  with $m$ edges and not containing a clique of order $t-1$  such that $\lambda(G)=\lambda_{(m, t-1)}^{3-}$.
\end{lemma}
{\em Proof of Theorem \ref{theorem 1}.} Let ${t-1 \choose 3} \le m \le {t-1 \choose 3} + {t-2 \choose 2}-\frac{1}{2}(t-1)$. Let $G$ be a 3-graph  with $m$ edges without containing $[t-1]^{(3)}$ such that $\lambda(G)=\lambda_{(m,t-1)}^{3-}$.  To prove Theorem \ref{theorem 1}, we only need to prove $\lambda_{(m,t-1)}^{3-}=\lambda(G)<\lambda([t-1]^{(3)}).$

By Lemma \ref{lemma1}, we can assume that $G$ is left-compressed. Let  $\vec{x}=(x_{1},x_{2},\ldots ,x_{n})$ be an optimal weighting for $G$. By Remark \ref{r1}(a), $x_1 \ge x_2 \ge \ldots \ge x_k >x_{k+1}=\ldots=x_{n}=0$. If $k\leq t-1$, $\lambda(G)<\lambda([t-1]^{(3)})$ since $G$ does not contain a clique order of $[t-1]$.  So we assume $k\geq t$. First we show that $k=t.$ Wee need the following lemma.
\begin{lemma} \cite{T} \label{LemmaTal9}
Let $G:=(V,E)$ be a left-compressed $3$-graph with $m$ edges such that $\lambda(G)=\lambda_{m}^{3}$. Let $b:=|E_{(k-1)k}|$. Let $\vec{x}:=(x_{1},x_{2},\ldots ,x_{k})$  be an optimal weighting for $G$ satisfying $x_1 \ge x_2 \ge \ldots \ge x_k >x_{k+1}=\ldots=x_{n}=0$. Then
$$|[k-1]^{(3)}\backslash E|\leq\lceil b(1+\frac{k-(b+2)}{k-3})\rceil.$$
\end{lemma}

Since $G$ is left-compressed and $1(k-1)k\in E$, then $\vert[k-2]^{(2)}\cap E_{k}\vert\geq 1$.
If $k\ge t+1$, then applying Lemma \ref{LemmaTal9},  we have $|[k-1]^{(3)}\backslash E|\leq k-2.$ Hence
\begin{eqnarray}\label{eqs31}
m=\vert E\vert &=&\vert E\cap [k-1]^{(3)}\vert +\vert [k-2]^{(2)}\cap E_k \vert +\vert E_{(k-1)k}\vert \nonumber\\
&\ge & {t \choose 3}-(t-1) +2 \nonumber\\
&\ge& {t-1 \choose 3} + {t-2 \choose 2}+1,
\end{eqnarray}
which contradicts to the assumption that $m\le {t-1 \choose 3} + {t-2 \choose 2}$. Recall that $k\ge t$, so we have
$$k=t.$$
Hence we can assume $G$ is on vertex set $[t]$.

Next we prove an inequality.
\begin{lemma}\label{Lemmaeq}  Let $G$ be a $3$-graph on the vertex set $[t]$. Let  $\vec{x}:=(x_{1},x_{2},\ldots ,x_{t})$ be an optimal weighting for $G$ satisfying $x_1 \ge x_2 \ge \ldots \ge x_t \ge 0$. Then $$x_1< x_{t-3}+x_{t-2}\  \rm or \ \lambda(G)\le \frac{1}{6}\frac{(t-3)^2}{(t-2)(t-1)}<\lambda([t-1]^{(3)}).$$
\end{lemma}
{\em Proof} If $x_1\geq x_{t-3}+x_{t-2}$,  then $$3x_1+x_2+\cdots+x_{t-4}>x_1+x_2+\cdots +x_{t-4}+x_{t-3}+x_{t-2}+x_{t-1}+x_{t}=1.$$ Recall that $x_1 \ge x_2 \ge \ldots \ge x_{t-4}$, we have $x_1> \frac{1}{t-2}$.
Using Lemma \ref{LemmaTal5}, we have
\begin{eqnarray*}
\lambda(G)&=&\frac{1}{3}\lambda(E_1,x)\leq\frac{1}{3}{t-1 \choose 2}(\frac{1-\frac{1}{t-2}}{t-1})^2\nonumber\\
&=&\frac{1}{6}\frac{(t-3)^2}{(t-2)(t-1)}
< \frac{1}{6}\frac{(t-3)(t-2)}{(t-1)^2}
=\lambda([t-1]^{(3)}).
\end{eqnarray*}
The first inequality follows from Theorem \ref{MStheo}.  Hence $\lambda(G)< \lambda([t-1]^{(3)})$, which contradicts to $\lambda(G)\geq\lambda([t-1]^{(3)})$. This completes the proof.\qed

The following lemma is proved in \cite{TPZZ2}.
\begin{lemma}({\cite{TPZZ2}, Lemma5.3})\label{Lemmatpzz} Let $G$ be a left-compressed 3-graph on the vertex set $[t]$. Let  $\vec{x}:=(x_{1},x_{2},\ldots ,x_{t})$ be an optimal weighting for $G$. Then $|[t-1]^{(3)}\backslash E|\leq t-3$, or  $\lambda(G)\leq \lambda ([t-1]^{(3)}).$
\end{lemma}
\begin{remark}\label{Remtpzz}
We can prove that $|[t-1]^{(3)}\backslash E|\leq t-3$, or  $\lambda(G)<\lambda ([t-1]^{(3)})$ under the condition of Lemma \ref{Lemmatpzz} through the method in \cite{TPZZ2}.
\end{remark}

Now we continue the proof of Theorem \ref{theorem 1}. Let $D=[t-1]^{(3)}\backslash E. $ Let $b=|E_{(t-1)t}|$.  By Lemma 
\ref{LemmaTal9}, we have $|D|\leq 2b.$ So $\lfloor \frac{|D|}{2}\rfloor\leq b$ and the triples $1(t-1)t,\cdots \lfloor\frac{|D|}{2}\rfloor(t-1)t$ are in $G$.  Let $G'=G\bigcup D\backslash \{1(t-1)t,\cdots \lfloor\frac{|D|}{2}\rfloor(t-1)t\}$. If $\lambda(G)<\lambda ([t-1]^{(3)}),$ we are done. Otherwise by Remark \ref{Remtpzz} we have $|D|\leq t-3$. So
\begin{eqnarray*}
|G'|&=&|G|+|D|-\lfloor\frac{|D|}{2}\rfloor\leq {t-1 \choose 3} + {t-2 \choose 2} - \frac{1}{2}(t-1)+t-3-\frac{t-3}{2}+1\nonumber\\
&=&{t-1 \choose 3} + {t-2 \choose 2}.
\end{eqnarray*}
 Note that $G'$ contains $[t-1]^{(3)}$. By Theorem \ref{TheoremPZ18}, we have $\lambda(G',\vec{x})\leq \lambda(G')=\lambda([t-1]^{(3)}).$

 Next we show that $\lambda(G,\vec{x})<\lambda(G',\vec{x}).$ By Remark \ref{r1}(b), $x_1=x_2=\cdots=x_{\lfloor\frac{|D|}{2}\rfloor}$. Hence
\begin{eqnarray*}
\lambda(G', \vec{x})-\lambda(G, \vec{x})&=&\lambda(D,\vec{x})-\lfloor\frac{|D|}{2}\rfloor x_1x_{t-1}x_{t}\nonumber\\
&\geq& |D|x_{t-3}x_{t-2}x_{t-1} -\lfloor\frac{|D|}{2}\rfloor x_1x_{t-1}x_{t}\nonumber\\
&>& |D|x_{t-3}x_{t-2}x_{t-1} -\lfloor\frac{|D|}{2}\rfloor (x_{t-3}+x_{t-2})x_{t-1}x_{t}.\nonumber\\
\end{eqnarray*}
In the last step, we used Lemma \ref{Lemmaeq}.  Recall that $x_1 \ge x_2 \ge \ldots \ge x_t>0,$ we have  $$|D|x_{t-3}x_{t-2}x_{t-1} -\lfloor\frac{|D|}{2}\rfloor (x_{t-3}+x_{t-2})x_{t-1}x_{t}\geq|D|x_{t-3}x_{t-2}x_{t-1}-|D|x_{t-3}x_{t-1}x_t\geq0.$$
Hence $\lambda(G,\vec{x})<\lambda(G',\vec{x})\leq \lambda([t-1]^{(3)})=\lambda(C_{3,m})$.  This completes the proof of Theorem \ref{theorem 1}.\qed
\section{Connection between the clique number and the Lagrangians of some left-compressed 3-graphs}
\label{sec:3}

In this section, we will  confirm Conjecture \ref{conjecture} and Conjecture \ref{conjecture2} for some left-compressed $3$-graphs with specified structures. As an application, we also obtain two weaker versions of Tur\'{a}n type result for left-compressed $3$-graphs.

\begin{theorem}\label{theorem 2}   Let $G:=(V,E)$ be a left-compressed $3$-graph on vertex set $[t]$ and $G$  does not contain a clique order of $\lfloor \frac{t-2}{2}\rfloor$. Then $$\lambda (G)\le \frac{1}{6}\frac{(t-3)^2}{(t-2)(t-1)}<\lambda ([t-1]^{(3)}).$$
\end{theorem}
\noindent{\em  Proof} 
The idea to prove Theorem \ref{theorem 2} is similar to that in  the proof of Lemma \ref{Lemmaeq}. Let $G:=(V,E)$ be a left-compressed 3-graph  with $m$ edges and $\omega(G)\leq \lfloor \frac{t-2}{2}\rfloor$. Recall $\omega(G)$ is the clique number of $G$. If $t\leq 5$, Theorem \ref{theorem 2} clearly holds. Next we assume $t\geq 6$. Let  $\vec{x}=(x_{1},x_{2},\ldots ,x_{t})$ be an optimal weighting for $G$ satisfying, $x_1 \ge x_2 \ge \ldots \geq x_{t}$. The clique number of $E_{t-3}$ must be smaller than $\frac{t-2}{2}$, otherwise $\omega(G)>\lfloor \frac{t-2}{2}\rfloor$ since $G$ is left-compressed. By Lemma \ref{Lemmaeq}, if $\lambda (G)>\frac{1}{6}\frac{(t-3)^2}{(t-2)(t-1)}$, we have $x_{t-3}> \frac{1}{2t}$. Using Lemma \ref{LemmaTal5} and Theorem \ref{MStheo}, we have
\begin{eqnarray*}
\lambda(G)&=&\frac{1}{3}\lambda(E_t,x)<\frac{1}{3}{\lfloor\frac{t-2}{2}\rfloor \choose 2}(\frac{1-\frac{1}{2t}}{\lfloor\frac{t-2}{2}\rfloor})^2\nonumber\\
&\leq&\frac{1}{6}\frac{t-4}{t-2}\frac{(2t-1)^2}{4t^2} \\
&< &\frac{1}{6}\frac{(t-3)^2}{(t-2)(t-1)}.
\end{eqnarray*}
which is a contradiction. This completes the proof.\qed

\begin{corollary}\label{coroturan1}   Let $G:=(V,E)$ be a left-compressed $3$-graph with $t$ vertices and $m$ edges. If $m\geq\frac{(t-3)^2t^3}{6(t-2)(t-1)}$,  then $G$ contains a clique order of $\lfloor \frac{t-2}{2}\rfloor$.
\end{corollary}
\noindent{\em  Proof} Let $G:=(V,E)$ be a $3$-graph with $t$ vertices and $m$ edges. Assume that $m\geq\frac{(t-3)^2}{6(t-2)(t-1)}t^3$. Clearly, $x_1=x_2=\cdots=x_t=\frac{1}{t}$ is a legal weighting for $G$. Hence $\lambda(G)\geq \frac{(t-3)^2}{6(t-2)(t-1)}t^3\frac{1}{t^3}=\frac{(t-3)^2}{6(t-2)(t-1)}$. However by  Theorem \ref{theorem 2} we know that $\lambda(G)<\frac{(t-3)^2}{6(t-2)(t-1)}$ if $G$ does not contain a clique order of $\lfloor \frac{t-2}{2}\rfloor$. This completes the proof. \qed


For the case of forbiding a clique of order $4$, we have the following result.

\begin{proposition} \label{proposition}
Let  $G$ be  a left-compressed 3-uniform graph on $[t]$ with $m$ edges.  If $G$ does not contain a clique of order $4$, then
$m  \le {2 \over 27}t^3.$
\end{proposition}

\noindent{\em Proof}
Let  ${\vec x}:=(x_1, x_2, \ldots, x_t)$ be an optimal vector of $G$.
We claim that all edges in $G$ must contain vertex 1. Otherwise, 234 is an edge of $G$ and  $G$ contains the clique $[4]^{(3)}$ since $G$ is left-compressed.   So
$$\lambda(G)\le x_1\cdot {1 \over 2}(x_2+x_3+\cdots+x_k)^2={1 \over 2}x_1(1-x_1)^2\le {1 \over 2} \times {4 \over 27}(x_1+{1-x_1 \over 2}+{1-x_1 \over 2})^3={2 \over 27}.$$

Let ${\vec y}=(y_1, y_2, \ldots, y_t)$ given by $y_i={1 \over t}$ for each $i, 1\le i\le t$. Then
${2 \over 27}\ge \lambda(G)\ge \lambda(G, {\vec y})={m \over t^3}.$ Therefore, $m  \le {2 \over 27}t^3.$ \qed


In \cite{BP2},  Bul\'o and Pelillo proved the following theorem.
\begin{theorem}(\cite{BP2})\label{theoBP}
An $r$-graph $G=(V,E)$ with $m$ edges and $t$ vertices, which contains no $p$-clique with $p\geq r$, then
$$m\leq{t \choose r}-\frac{t}{(r-1)r}[(\frac{t}{p-1})^{r-1}-1].$$
\end{theorem}

\begin{remark}  (1) We note that  Theorem \ref{theorem 2} and Corollary \ref{coroturan1} establish a connection between Lagrangian and clique number for $3$-graphs. They also provide evidence for Conjecture \ref{conjecture2}.

(2) For the case $r=3$ and $p=\lfloor\frac{t-2}{2}\rfloor$, the upper bound in Theorem \ref{theoBP} is bigger than $\frac{(t^4-11t^3+39t^2-72t+48)t}{6(t-4)^2}$. Since $$\frac{(t^4-11t^3+39t^2-72t+48)t}{6(t-4)^2}>\frac{(t-3)^2t^3}{6(t-2)(t-1)}$$ when $t\geq 38$, the result in Corollary \ref{coroturan1} is better than the result in Theorem \ref{theoBP} under the left-compressed condition.

(3) Again, for the case $r=3$ and $p=4$, the upper bound in Theorem \ref{theoBP} is bigger than the bound in Propostion \ref{proposition} under the left-compressed condition.
\end{remark}

Next we give the following partial result to Conjecture \ref{conjecture}.

\begin{theorem}\label{theorem43} Let $m$, $t$, and $a$ be positive integers satisfying $m={t-1 \choose 3} + {t-2 \choose 2}+a$ where $1\leq a\leq t-2$.  Let $G=(V,E)$ be a left-compressed 3-graph on the vertex set [t] with $m$ edges satisfying  $|E_{(t-1)t}|\le \frac{2t+3a-4}{5}$. If $G$ contains a clique of order $t-1$, then $\lambda (G)\leq\lambda (C_{3,m})$.
\end{theorem}

\noindent{\em Proof} Let $G$ be a $3$-graph with $m$ edges and  containing a clique of order  $t-1$.
Assume $\vec{x}:=(x_{1},x_{2},\ldots ,x_{t})$  is an optimal weighting for $G$ satisfying $x_1 \ge x_2 \ge \ldots \ge x_t \ge 0$. We will prove that  $\lambda(C_{3,m},\vec{x})-\lambda(G,\vec{x})\geq 0$. Therefore $\lambda(C_{3,m})\geq \lambda(C_{3,m},\vec{x})\geq \lambda(G, \vec{x})=\lambda(G)$.

By Remark \ref{r1}(b) and noting that $G$ contains $[t-1]^{(3)}$, we have
\begin{eqnarray}\label{eqa}
x_{1}&=& x_{t-3}+\frac{\lambda(E_{1\setminus(t-3)},\vec{x})}{\lambda(E_{1(t-3)},\vec{x})}
=x_{t-3}+\frac{\lambda(E_{t-3}^c,\vec{x})}{\lambda(E_{1(t-3)},\vec{x})}
=x_{t-3}+\frac{x_t\lambda(E_{(t-3)t}^c,\vec{x})}{\lambda(E_{1(t-3)},\vec{x})},
\end{eqnarray}
and
\begin{eqnarray}\label{eqb}
x_{t-2}&=& x_{t-1}+\frac{\lambda(E_{(t-2)\setminus(t-1)},\vec{x})}{\lambda(E_{(t-2)(t-1)},\vec{x})}
=x_{t-1}+\frac{x_t\lambda(E_{(t-2)t}\bigcap E_{(t-1)t}^c,\vec{x})}{\lambda(E_{(t-2)(t-1)},\vec{x})}.
\end{eqnarray}

Let $b:=|E_{(t-1)t}|$. Since $G$ contains the clique $[t-1]^{(3)}$, we have $|[t-2]^{(2)}\backslash E_t|=b-a.$ Note that $|E_{(t-3)t}^c|\leq |E_{(t-2)t}^c|$ since $G$ is left-compressed, we have $|E_{(t-3)t}^c|\leq \frac{b-a}{2}+1$(Note that $t-1\in E_{(t-3)t}^c$). On the other hand, $|E_{(t-2)t}|=t-2-|E_{(t-2)t}^c|\geq t-2-(b-a)-1$(Note that $t-1\in E_{(t-2)t}^c$) and $|E_{(t-2)t}\bigcap E_{(t-1)t}^c|\geq (t-2)-(b-a)-1-b=t-2b+a-1$. Recalling that $|E_{(t-1)t}|\le \frac{2t+3a-4}{5}$, we have  $|E_{(t-3)t}^c|\leq |E_{(t-2)t}\bigcap E_{(t-1)t}^c|$. Let $i$ be the minimum integer in  $E_{(t-3)t}^c$ and $j$ be the minimum integer in  $E_{(t-2)t}\bigcap E_{(t-1)t}^c$. Because $G$ is left-compressed, we have $i\geq j.$ Hence
\begin{eqnarray}\label{eqc}
\lambda(E_{(t-3)t}^c,\vec{x})\leq \lambda(E_{(t-2)t}\bigcap E_{(t-1)t}^c,\vec{x}).
\end{eqnarray}
Since $x_1 \ge x_2 \ge \ldots \ge x_t$.
Next we prove that
\begin{eqnarray}\label{ineqadd}
\lambda(E_{1(t-3)},\vec{x})-\lambda(E_{(t-2)(t-1)},\vec{x})&=& x_{t-2}+x_{t-1}+x_{t}-x_{1}-x_{t-3} \geq 0.
\end{eqnarray}
To verify (\ref{ineqadd}), by Remark \ref{r1}(b), we have
\begin{eqnarray}\label{By1}
x_{1}=x_{t-1}+\frac{\lambda(E_{1\backslash (t-1)},\vec{x})}{\lambda(E_{1(t-1)},\vec{x})}\leq x_{t-1}+\frac{(x_{2}+\cdots+x_{t-2})x_{t}}{x_{2}+\cdots+x_{t-2}+x_{t}}\leq x_{t-1}+x_{t};
\end{eqnarray}

 \begin{eqnarray}\label{By2}x_{1}&=&x_{t-2}+\frac{\lambda(E_{1\backslash (t-2)},\vec{x})}{\lambda(E_{1(t-2)},\vec{x})}\nonumber \\
 &=&x_{t-2}+\frac{\lambda(E_{(t-2)t}^c,\vec{x})}{1-x_{1}-x_{t-2}}x_{t}\nonumber \\
 &\leq&x_{t-2}+\frac{\lambda(E_{(t-2)t}^c,\vec{x})}{1-x_{t-3}-x_{t-1}-x_{t}}x_{t};
\end{eqnarray}
and
\begin{eqnarray}\label{By3}x_{t-3}&=&x_{t-1}+\frac{\lambda(E_{(t-3)\backslash (t-1)},\vec{x})}{\lambda(E_{(t-3)(t-1)},\vec{x})}\nonumber \\
&=&x_{t-1}+\frac{\lambda(E_{(t-3)t}\bigcap E_{(t-1)t}^c,\vec{x})}{1-x_{t-3}-x_{t-1}-x_{t}}x_{t}.
\end{eqnarray}
Adding (\ref{By2}) and (\ref{By3}), we obtain that
\begin{eqnarray*}
x_{1}+x_{t-3}&\leq& x_{t-2}+ x_{t-1}+\frac{\lambda(E_{(t-2)t}^c,\vec{x})+\lambda(E_{(t-3)t}\bigcap E_{(t-1)t}^c,\vec{x})}{1-x_{t-3}-x_{t-1}-x_{t}}x_{t}.
\end{eqnarray*}
Clearly $t-3 \notin E_{(t-3)t}$. Since $G$ is left-compressed and $G\neq C_{3,m}$,  we have $t-2 \notin E_{(t-3)t}$. On the other hand both $t-3$ and $t-2$ are in $E_{(t-1)t}^c$. Hence $|E_{(t-2)t}^c|+|E_{(t-3)t}\bigcap E_{(t-1)t}^c|\leq |E_{(t-2)t}^c|+|E_{(t-1)t}^c|-2$. Recalling that $|E^c|\leq t-3$, we have $|E_{(t-2)t}^c|+|E_{(t-1)t}^c|\leq|E^c|\leq t-2$(Note that $t-1\in E_{(t-2)t}^c$ and $t-2\in E_{(t-1)t}^c$ ) and $|E_{(t-2)t}^c|+|E_{(t-3)t}\bigcap E_{(t-1)t}^c|\leq t-4$.  Clearly  $b\geq 2$.
Hence $2$ is not in $E_{(t-1)t}^c$ and $E_{(t-2)t}^c$. Recalling that $x_1 \ge x_2 \ge \ldots \ge x_t$, we have $\frac{\lambda(E_{(t-2)t}^c,\vec{x})+\lambda(E_{(t-3)t}\bigcap E_{(t-1)t}^c,\vec{x})}{1-x_{t-3}-x_{t-1}-x_{t}}\leq 1$. So, (\ref{ineqadd}) is true.
This implies that $\lambda(E_{(t-2)(t-1)},\vec{x})\leq\lambda(E_{1(t-3)},\vec{x})$.
Combining (\ref{eqa}),\ref{eqb}) and (\ref{eqc}), we obtain that
$x_1-x_{t-3}\leq x_{t-2}-x_{t-1}$ and $x_{t-3}x_{t-2}x_t-x_1x_{t-1}x_{t}\geq 0$.
Hence
\begin{eqnarray}\label{eq37}
\lambda(C_{3,m},\vec{x})-\lambda(G,\vec{x})&=&\lambda([t-2]^{(2)}\backslash E_t,\vec{x})-|[t-2]^{(2)}\backslash E_t|x_{1}x_{t-1}x_{t}\nonumber \\
&\geq& |[t-2]^{(2)}\backslash E_t|(x_{t-3}x_{t-2}x_t-x_1x_{t-1}x_t)\nonumber \\
&\geq& 0.
\end{eqnarray}
 This completes the proof. \qed

\begin{acknowledgements}
 This research is partially  supported by National Natural Science Foundation of China (No. 11271116).
\end{acknowledgements}



\end{document}